\documentclass[11pt]{article}
\usepackage[utf8]{inputenc}

\usepackage[final]{graphicx}

\usepackage[a4paper, textwidth={14cm}]{geometry}

\usepackage{amsthm}
\usepackage{amsmath}
\usepackage{amsfonts}
\usepackage{amssymb}

\newcommand{\highlightcol}{blue}
\usepackage{hyperref}
\hypersetup{
  colorlinks   = true,    
  urlcolor     = \highlightcol,    
  linkcolor    = \highlightcol,    
  citecolor    = \highlightcol     
}

\usepackage[font=small,labelfont=bf]{caption}
\usepackage{subcaption}

\newtheorem*{theorem*}{Theorem}
\newtheorem*{definition*}{Definition}

\newcommand{\keywords}[1]{{\bf Keywords:} #1}

\newcommand{\remoteReference}[1]{#1}

\begin{document}

\title{Formalising the intentional stance 1:\\ attributing goals and beliefs to stochastic processes}

\author{Simon McGregor${}^1$, timorl, Nathaniel Virgo${}^{2,3}$ \\[1ex]
\small
${}^1\,$University of Sussex, UK \\
\small
${}^2\,$Centre of Data Innovation Research,\\
\small School of Physics, Engineering \& Computer Science,\\ 
\small University of Hertfordshire, UK\\
\small
${}^3\,$Earth-Life Science Institute,\\ \small Institute of Science Tokyo, Japan.
}

\maketitle

\begin{abstract}
	This article presents a formalism inspired by Dennett's notion of the intentional stance. Whereas Dennett's treatment of these concepts is informal, we aim to provide a more formal analogue. We introduce a framework based on stochastic processes with inputs and outputs, in which we can talk precisely about
\emph{interpreting} systems as having \emph{normative-epistemic states}, which combine belief-like and desire-like features. Our framework is based on optimality but nevertheless allows us to model some forms of bounded cognition.

One might expect that the systems that can be described in normative-epistemic terms would be some special subset of all systems, but we show that this is not the case: every system admits a (possibly trivial) normative-epistemic interpretation, and those that can be \emph{uniquely specified} by a normative-epistemic description are exactly the deterministic ones.
Finally, we show that there is a suitable notion of Bayesian updating for normative-epistemic states, which we call \emph{value-laden filtering}, since it involves both normative and epistemic elements. For unbounded cognition it is always permissible to attribute beliefs that update in this way. This is not always the case for bounded cognition, but we give a sufficient condition under which it is.

This paper gives an overview of our framework aimed at cognitive scientists, with a formal mathematical treatment given in a companion paper.

\end{abstract}

\newcommand{\inputSpace}{I}
\newcommand{\outputSpace}{O}
\newcommand{\transducer}{\pi}
\newcommand{\elementString}[1]{\mathbf{#1}}
\newcommand{\stateSpace}{S}
\newcommand{\actionSpace}{A}
\newcommand{\noSuccess}{\bot}
\newcommand{\success}{\top}
\newcommand{\environment}{\varepsilon}
\newcommand{\policy}{\transducer}


\section{Introduction}

We are interested in the question of \emph{what does it mean for a system to be an agent?}
More specifically, we are interested in the relationship between what we might call \emph{mechanistic descriptions} of a system (explanations of its behaviour in terms of its internal state) and what we call \emph{normative-epistemic descriptions}.
These are descriptions in terms of a combination of norms and beliefs.
For example, we might ascribe a goal to the system and beliefs about its surroundings, and explain its behaviour in terms of those.
Both types of description concern the externally observable behaviour of a system, but they are quite different in nature, and we want to understand where this difference comes from and exactly how they relate.
We take a position similar to \cite{Dennett1981Book} in that we regard in intentional descriptions as a \emph{stance}, i.e.\ an optional perspective that can be taken by an observer, rather than an inherent property of a system.
Our goal is to express this idea mathematically.

In order to do this, we propose a simple mathematical framework in which systems can be described in both mechanistic and normative-epistemic terms, and the relationships between the two can be understood.
We don't aim for the most general framework possible but for the simplest one in which our main points can be expressed.
We describe the framework intuitively in this paper, with the technical details present in the companion work.

The current work follows a particular conceptual methodology, closely related to Dennett's intentional stance; we describe this methodology, and its relations to Dennett's work in Section~\ref{sec.as-if}. Section~\ref{sec.results} explains the mathematical framework we are proposing, and summarises the mathematical results described in the companion paper.
In Section~\ref{sec.literature} we describe relations to previous formalisms: namely, representation theorems for expected utility theory, the free energy principle, computational mechanics and inverse reinforcement learning. Finally, we discuss limitations of our framework and possible future work in Section~\ref{sec.limitations}.

\section{The `As-If' Approach And Dennett's Intentional Stance}
\label{sec.as-if}

The current work is part of a research programme of `as-if' agency.
The idea of the `as-if' approach \cite{McGregor2016, McGregor2017} is to formally address what it means to treat a system as if it were an agent, while remaining agnostic about whether or not there are additional criteria that must be met for a system to be deemed a `true' agent.

This idea can be seen as an attempt to formalise Dennett's intentional stance \cite{Dennett1975, Dennett1981, Dennett2006}; in doing so, we are also able to formally model relationships between (analogues of) Dennett's intentional stance and his physical stance.
The physical stance relates to technical concepts of stochastic Moore machines and unifilar machines, described in the companion paper, while the intentional stance corresponds to our exploration of the relationship between transducers and teleo-environment, in Section~\ref{sec.teleo} and throughout both papers.

We emphasise, though, that while our mathematical treatment has some
features that are strongly suggestive of Dennett’s ideas, and is inspired by Dennett’s work, we don’t consider it a full formalisation of Dennett’s work as such, and it
isn't intended to address
everything Dennett is concerned with.

In the past there have been significant disagreement about what sorts of system should be regarded as `agents'. Some authors (e.g.\ \cite[p.34]{Russel2009} or \cite{Wooldridge1994}) categorise even simple feedback systems (such as thermostats or guided torpedoes) as agents. Others (e.g.\ \cite{Dretske1999} or \cite{Jonas1953}) explicitly reject such a characterisation, claiming that thermostats may behave (in some ways) as if they were agents, but are not truly agents. A discussion of this general debate can be found in \cite{Martinelli2023}. Currently, echoes of this disagreement can be seen in disputes about the extent to which AI systems such as large language models (LLMs) should be understood as agents.
The idea of `as-if' agency is that we should not ignore simple systems such as thermostats when studying agency, even if the skeptics are correct and thermostats lack some important property. It is still valuable to have a theory that explains in what sense even a thermostat behaves \textit{as though} it were an agent.

The `as-if' approach involves considering the match between a system's physical behaviour, and the prescriptions of a normative theory of `ideal' rationality \cite{McGregor2016,McGregor2017}.
As outlined in \cite{McGregor2017},
``[I]f the causal effects exerted by some system $X$ on its environment are roughly rational decisions under some formal theory of rationality $R$ parameterised by some cognitive parameters $C$, we will say that $X$ is an as-if agent with respect to $(R, C)$.''

We will not define what a `formal theory of rationality' is, but our notion of
teleo-environment (Section \ref{sec.teleo}) will serve as an example.

The theory we describe in this paper is not utility theory (as suggested in \cite{McGregor2017}), but instead a simplified normative model of goal-belief rationality, in which the agent is treated as if it is trying to maximise the probability with which a goal event will occur.
Some philosophical motivation for this is given in \cite{McGregor2020mandate}.

More recently, \cite{VirgoBiehlMcGregor21} introduced a notion of `consistent Bayesian interpretation.'
This is a slightly different version of the `as-if' approach, in which the agent's internal states are interpreted as having a semantic meaning, rather than only its externally observable behaviour.
It introduces an important element, which is that it makes sense to ask for an `as-if' interpretation to be consistent over time.
In particular, if we can treat a system as if it has a Bayesian prior at time $t$, then at time $t+1$, after having received some new sensory information, we should be able to treat the system as having a new set of beliefs, given by the appropriate Bayesian posterior.
Ref.\ \cite{VirgoBiehlMcGregor21} shows that the conditions for such a `consistent Bayesian interpretation' to be possible are relatively weak, at least in formal terms.
When such interpretations exist they are not unique, which is an important aspect of the `as-if' approach in general: there is rarely only one way in which a system can reasonably be treated as if it were an agent.
This approach is extended to rational behaviour using value functions in \cite{Biehl2022}.

In the current paper and its companion, we are also concerned with consistency.
We show that we obtain a form of consistency (Section \ref{sec.filtering}) if
certain assumptions are met. The form of consistency is not what we
might naively expect, since the update rule is given by `value-laden filtering'
(Theorem \ref{thm.filtering}), rather than vanilla Bayesian filtering.
We also observe that this `value-laden filtering' consistency does not hold in general in bounded-rational contexts.
In particular, in Section \ref{sec.bounded}, we describe a case in which it doesn't hold, due to constraints on memory, and mention another consequence of optimality under memory constraints, namely that the uniquely optimal policy may be non-deterministic, as previously observed in the absent-minded driver problem (\cite{PICCIONE19973}) and discussed more generally in \cite{Icard2021}.

\section{The Teleo-Environment Framework}
\label{sec.results}

In this section we describe our mathematical framework intuitively, and provide a detailed summary of our conceptually relevant results, using a minimum of technical language. We encourage the
reader interested in technical details to read the companion paper, which provides a rigorous treatment.

Consider the following situation:
a person who we will refer to as `the theorist' is studying a physical system.
Although we use the term `physical system'\!, a suitable image to have in mind might be that of an animal, or even better, a robot.
(We can imagine that the robot was designed by someone other than the theorist, for a purpose that the theorist doesn't necessarily know about.)
The theorist's goal is to understand whether the system they are dealing with is an agent.
Or, more accurately, whether it can be consistently interpreted \emph{as if} it were an agent, in the sense that it takes actions that are consistent with trying to achieve some particular goal.
There might be more than one way in which the system can be so interpreted; if so, the theorist is interested in all of them.

\subsection{Transducers}
In order to reason about physical systems in such a situation, we make the following assumptions:
the system has a set of possible outputs or `actions' $\actionSpace$ and a set of possible inputs or `sensor values' $\stateSpace$ (both finite).
Time proceeds in discrete steps and the system's behaviour may be stochastic.
Classically, we would model its externally observable behaviour in terms of a \emph{controlled stochastic process}, that is, a conditional probability distribution over infinite sequences of outputs or \emph{actions} $a_1, a_2, \dots$ given infinite sequences of inputs or \emph{sensor values} $s_1, s_2, \dots$.
This conditional distribution must have the property that the output at time $t$ is independent of all inputs received after time $t-1$.
(See \cite{Hutter2005universal} for a similar approach.)
We refer to this as the \emph{causality condition}.

However, we choose instead to model it using an object we call a \emph{transducer}, which could be
described as follows.
\begin{definition*}[\textbf{\remoteReference{3.0.1}} in the companion paper]
  Given finite sets of inputs $\inputSpace$ and outputs $\outputSpace$ a transducer (from $\inputSpace$ to $\outputSpace$) is a
  mathematical object that provides a distribution over $\outputSpace$ and
  deterministically changes into another transducer when provided with a
  sensory value and an action chosen from its distribution.
\end{definition*}

If we start with a transducer $\pi$, give it an input $i$ and observe that it produced the output $o$, we obtain a new transducer that we refer to as ``$\pi$ evolved by $(i,o)$.''
Intuitively a transducer represents a process that simultaneously receives some input
data and stochastically produces an output, after which it is ready for another
step, perhaps in a different state.

This definition may seem not be well founded, since we are defining transducers in terms of themselves without a base case.
However, a mathematical framework known as coinduction allows such definitions to be treated rigorously.
The companion paper uses coinduction to reason about transducers defined in this way, and briefly explores the relation between transducers and controlled stochastic processes.

We assume that the theorist has access to a \emph{policy}  representing their system's observable behaviour.
This is simply a transducer
from $\stateSpace$ to $\actionSpace$.
This means that the theorist is not limited to observing a single sample of the system's behaviour given a single sequence of inputs; instead they know\footnote{or at least, have a model of}
how the system \emph{would} behave
given any number of different input sequences, over any number of trials.

\subsection{Teleo-environments}
\label{sec.teleo}

\begin{figure}[t]
  \centering
  \begin{subfigure}[b]{0.45\textwidth}
    \centering
    \includegraphics[width=\textwidth]{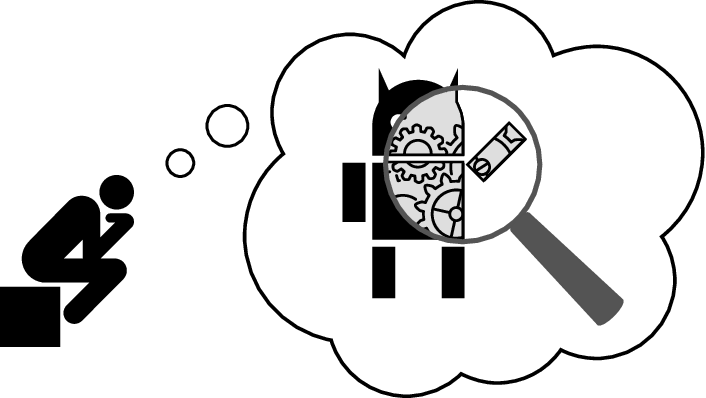}
    \caption{mechanistic description}
    \label{mechanistic.fig}
  \end{subfigure}
  \hfill
  \begin{subfigure}[b]{0.45\textwidth}
    \centering
    \includegraphics[width=\textwidth]{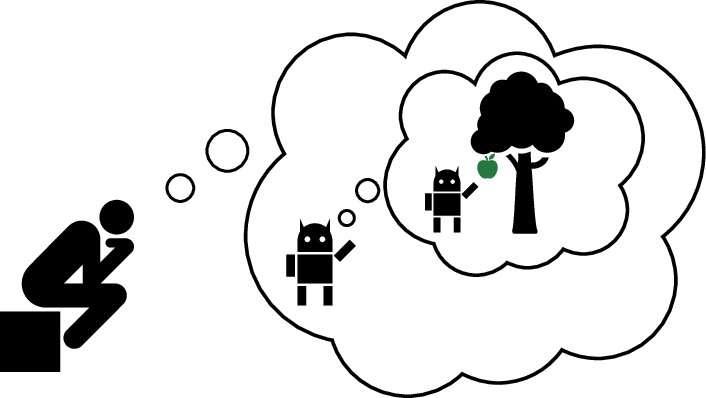}
    \caption{normative-epistemic description}
    \label{normative-epistemic.fig}
  \end{subfigure}
  \caption{Two ways in which a theorist might describe the behaviour of a
    physical system. $(a)$ in a mechanistic description the theorist posits an
    internal state to the system such that its dynamics can explain the
    observable behaviour; in the case of a robot this could include hypothesised
    workings of internal control systems in addition to the literal mechanical
    gears depicted here. We formalise mechanistic descriptions in terms of
    stochastic Moore machines in the companion paper.
    $(b)$ in a normative-epistemic description the theorist posits an
    environment that the system believes itself to be in (according to the
    theorist's interpretation) and a goal that it seeks to achieve --- in this
    case the apple --- such that the system's behaviour can be explained by
    acting optimally in order to reach the goal. Normative-epistemic
    descriptions are the focus of our paper and are formalised in terms of
    teleo-environments.
    We emphasise that these are two different kinds of explanation for the same behaviour of the same system, and not two different types of system or two competing hypotheses.
  }
  \label{fig.robot}
\end{figure}

Given knowledge of a system's observable behaviour in the form of a policy, one question one could ask is, ``given the known behaviour, what kinds of physical mechanism might exist inside the system that would give rise to its behaviour?''
Although we briefly discuss this question in Section~\remoteReference{3.3.4} of the companion paper, it will not be our main focus.

Instead, we are concerned with the question ``is the system's behaviour consistent with taking optimal actions in order to achieve some  goal?''
The relationship between these two questions is illustrated in Figure~\ref{fig.robot}.

To address this, we introduce the notion of \emph{teleo-environment}.
This is also a transducer, but with the opposite interface to the system:
where the system takes input from a set $\stateSpace$ and gives outputs in a set $\actionSpace$, its environment takes actions $\actionSpace$ as inputs and gives sensor values $\stateSpace$ as outputs, so that the two systems can be coupled.
(They are assumed to couple in such a way that they give their respective outputs simultaneously, as opposed to alternating.)

In addition to this, a teleo-environment comes equipped with a notion of \emph{success}.
The idea is that success is an event that may or may not occur on any given time step, and the agent's goal (as attributed by the theorist) is to maximise the probability that success occurs \emph{at least once}.
The agent does not necessarily know whether success has been achieved or not on any given time step, and hence can be `motivated' to continue `trying' to succeed even if success has already occurred.
This differs somewhat from the usual formulation of partially observable Markov decision processes (POMDPs) in terms of exponentially discounted reward, but is arguably more intuitive and results in an arguably simpler mathematical framework; it also has some interesting consequences.
Hence we define a teleo-environment as a transducer from $\actionSpace$ to $\stateSpace \times \{ \noSuccess, \success \}$ where $\noSuccess$ is a `lack of success' signal, and $\success$ is a `success' signal.
This notional `telos' channel is not directly observable by the system, but forms part of our concept of optimality.

It is worth noting that our notion of `goal' is inherently tied to a notion of `environment', since teleo-environments embody both in a way that can't easily be separated.
When we vary the goal in order to ask about the set of all goals for which a given agent is optimal, we vary the environment as well.
(This is in contrast to the usual situation in inverse reinforcement learning, in which the environment and goal are separate, and the environment is known but the goal is not.)

Our position is that, intuitively, a teleo-environment can be taken to represent a particular normative-epistemic state, meaning a description of a system's behaviour in terms of goals and beliefs.
We regard it as describing
both \emph{what the agent believes} about its environment (as attributed to it by the theorist), and \emph{what it is (interpreted as) trying to do}.
The idea is that if a policy $\policy$ is optimal for a teleo-environment $\environment$, it's consistent to claim that $\policy$ behaves as it were an ideally rational agent, with beliefs and goal represented by $\environment$. In this sense, we can treat $\environment$ as a permissible `mental state attribution' for $\policy$.
The teleo-environment $\environment$ does not necessarily represent the dynamics of the true environment, but rather the environment that the transducer~$\policy$ `believes' it is in, according to the {normative-epistemic state} attributed to it by the theorist.

Since $\environment$ is not unique, multiple distinct mental state attributions may be applied to the same system.
We see this as a strength of our approach rather than a weakness, since it captures the idea that not only is taking the intentional stance a choice, but that there can be choices in how it is applied as well.

\subsection{Filtering and the Bellman property}
\label{sec.filtering}
Within this simplified model of ideal rationality, we can ask: as a system `evolves' its policy through interaction with its environment, how do the permissible mental state attributions change?

We show that, if $\policy$ is optimal for $\environment$, then after evolution with a sensorimotor pair $(s, a)$ the evolved policy $\policy'$ is optimal for a similar evolution $\environment'$ of $\environment$. We argue in that section that this can be seen as a form of \textit{Bayesian filtering}.
In other words, among the many beliefs that are permissibly attributable to the system after it receives new information, one is guaranteed to look like an update of the previous attributed belief according to Bayes' rule.
(There is some subtlety about exactly what data is conditioned on, however, which will be discussed shortly.)

We propose a general intuition along the following lines. Consider any notion of optimality for a control policy $x$ in an environment $y$, and some notion of updating over time for the policy and environment (we will be deliberately vague about what this `updating' means in general). Suppose something like Bellman's principle holds, namely that if a policy $x$ is optimal for $y$, then the evolution $x'$ by a single time step will be optimal for the evolution $y'$ of $y$. Again, $y$ will function as a permissible `mental state' attribution for $x$ (in the sense that no inconsistency arises from supposing that $x$ is an ideal agent trying to solve the problem represented by $y$) and similarly for $y'$ and $x'$.
Consequently, the `update' operator that we apply to $y$
to obtain $y'$ can be seen as updating beliefs attributed to $x$, to obtain beliefs that can be attributed to $x'$, given any information received by $x$ during the update.
In this sense the `update' rule for $y$ can be seen as a model of idealised rationality, which generalises Bayes' rule in some sense.

When we work through the details for our teleo-environment model, we obtain an interesting insight. One might reasonably guess that a Bayesian model of rational mental state update might emerge from our model. In fact it does, but with an important nuance: so to speak, the system is legitimated to update its `belief' in a Bayesian manner, conditioning on the sensorimotor signals observed in the previous time step, provided that it additionally conditions on not having achieved its goal in the previous time step.
That is, it conditions on its received sensor value $s\in S$, but it also conditions on $\noSuccess$, the lack-of-success signal, even though it doesn't actually have access to the success signal, and indeed the true value of the success signal might have been $\success$.
We call this `value-laden filtering', and it arises from the relevant
Bellman-like principle.

\begin{theorem*}[\textbf{Corollary \remoteReference{5.0.3}} in the companion paper]
\label{thm.filtering}
  If a policy $\policy$ is optimal for an environment $\environment$, then
  given a possible sequence of actions $\elementString{a}$ and sensor values
  $\elementString{s}$ the policy $\policy'$, is optimal for the environment
  $\environment'$.
  Here, $\policy'$ is $\policy$ evolved by
  $\elementString{a}$ and $\elementString{s}$, and
  $\environment'$ is  $\environment$ evolved by $\elementString{a}$ and $\elementString{s}$
  (combined with a sequence of the lack-of-success symbol $\noSuccess$).
\end{theorem*}

This is a little surprising, since an external theorist, observing the sensorimotor interactions between the system and its environment (without any additional information about whether a goal had been achieved), and applying Bayesian reasoning themselves to predict future interactions, would condition only on the sensorimotor signals.
In this way, the properties of the optimality criterion in our model mean that we can derive a permissible rule for pragmatically optimal `mental state' change (namely, value-laden filtering) which differs subtly from the standard Bayesian picture, whilst still being recognisable as Bayesian filtering.

Value-laden filtering conditions on the lack-of-success signal in addition to the received signal.
Corollary \remoteReference{5.0.3} says that if a normative-epistemic state is a permissible attribution on one time step, then an updated version of the same normative-epistemic state will be a permissible attribution on the next time step, where the updating is given by value-laden filtering.

We are using slightly circumspect language here because of a subtle issue.
Because we use ``permissible attribution'' to mean any teleo-environment $\varepsilon$ for which a given $\pi$ is optimal, there are in general many permissible mental state attributions for a given policy $\pi$.
We might want to talk about an agent's `current beliefs' and then say that after the agent has undergone some sensor-motor interaction with its environment (i.e. been evolved by some $(\bf a,s)$), the resulting transducer $\pi'$ has an updated set of beliefs given by $\varepsilon'$.
However, this doesn't quite make sense, because in general there are many possible mental state attributions for $\pi'$, so ``its beliefs'' are not unique.
The above theorem says that the updated beliefs $\varepsilon'$ (obtained by value-laden filtering) are always \emph{among} the permissible mental state attributions for $\pi$, but there may be (and in general are) many others.

For \emph{some} permissible beliefs $\varepsilon$, it happens that the updated belief $\varepsilon'$ is the same as the one that would be obtained not by value-laden filtering but by ordinary Bayesian filtering. That is, one can obtain $\varepsilon'$ by evolving $\varepsilon$ by $\bf a$ and $\bf s$ only, without the additional sequence of lack-of-success symbols.

In fact, given any permissible attribution $\varepsilon$, there is a way to construct another attribution, which we write $Z(\varepsilon)$, such that (i) for every policy $\pi$, the probability of achieving success at least once is equal to that for $\varepsilon$, and (ii) $Z(\varepsilon)$ has the property that $Z(\varepsilon)$-evolved-by-$\bf s$-and-$\bf a$ is a permissible attribution to $\pi'$,
i.e.\ we can choose to regard $Z(\varepsilon)$ as updating by ordinary Bayesian filtering instead of by value-laden filtering.
The attribution $Z(\varepsilon)$ is defined in Definition \remoteReference{5.2.3} in the companion paper, and the proof of this statement is in Corollary \remoteReference{5.2.7}.

To give slightly more detail, $Z(\varepsilon)$ is obtained from $\epsilon$ by a process called ``single-success truncation.''
It is identical to $\varepsilon$ except success will only ever occur at most once.
If success does occur then $Z(\varepsilon)$ evolves into an environment that behaves identically to $\varepsilon$, except that the success signal does not occur when it otherwise would.

Both truncating environments, as well as value-laden filtering serve essentially
the same purpose in modifying the update mechanism of attributed states. A
transducer only represents the current state and the future of the system,
completely disregarding the past. However, if an agent already achieved success
in a past step, it should no longer care whether it will achieve success again.
For example, consider an environment in which, on the first time step, the agent either achieves success and enters environment $X$, or doesn't achieve success and enters environment $Y$.
Environment $X$ is such that the agent must take action $A$ in order to achieve success (for the second time), while environment $Y$ is such that the agent must take action $B$.
And suppose further that the agent's sensors do not allow it to distinguish between the two environments.
Even if environment $X$ is entered much more often, say with probability 0.9, an optimal
agent should still ignore it and always take action $B$.
The point is that in this situation, the optimal behaviour from the second time step onwards depends on the probability that success has already occurred, which is correlated with the probability of being in environment $X$ versus environment $Y$.

From the perspective of a person trying to solve this task, the reasoning at the second time step should be something along the lines of ``either I am in environment $X$, in which case I can achieve success by doing action $A$, but this doesn't matter, since I've already achieved success in that case, or I'm in environment $Y$, in which case I have not achieved success yet and must take action $B$ in order to achieve it at least once.
Since it's only in the case of environment $Y$ that my actions matter at all, I can discount the possibility being in environment $X$ entirely and behave as if I am in environment $Y$.''

If we were to model the agent's knowledge as a teleo-environment  but update it using vanilla Bayesian updating, it would result in an agent that always tries to achieve success starting from its current situation, regardless of whether success has already been achieved in a previous time step.
Truncation avoids this by letting the agent naturally ignore any futures where it would
have achieved success again, while value-laden filtering avoids it by explicitly excluding
the futures that follow a success from the model.

In future, it may be of interest to explore further how different pragmatic notions of optimality can induce different permissible update rules for mental state, and how these rules relate to `purely normative' models of ideal mental state update.

\subsection{Specification}

Another question one might want to ask is whether a given behaviour can be \emph{specified} by a given teleo-environment.
That is, for a system with given behaviour, does there exist a teleo-environment for which it is uniquely optimal?
Our results here are as one might expect from a framework based on utility-maximisation.
Firstly, a behaviour can only be uniquely optimal for a given teleo-environment if it is deterministic, since otherwise it would be a mixture of two or more deterministic behaviours, each of which would also be optimal.
Secondly, the specifiable behaviours are \emph{exactly} the deterministic ones: if the behaviour is deterministic then there exists a teleo-environment for which it is uniquely optimal.
(There may be many such teleo-environments, but in particular there is always one in which the system must perform exactly the specified behaviour, otherwise it immediately permanently loses all hope of achieving success.)

\begin{theorem*}[\textbf{\remoteReference{6.2.3}} and \textbf{\remoteReference{6.3.2}} in the companion paper]
\label{thm.deterministicAreSpecifiable}
  Any deterministic policy $\policy$ is specifiable, i.e. there exists a
  teleo-environment for which $\policy$ is uniquely optimal, and furthermore any specifiable policy is deterministic.
\end{theorem*}

Thus, at least when it comes to deterministic behaviour, there is no difference in the set of behaviours that can be specified in a normative-epistemic way and the set of behaviours that can be specified by a framework based on internal state.

This is in a sense a `no-go theorem'.
One might initially have the intuition that systems that are `agents' (in the sense of being optimal for some teleo-environment) might have special features to their behaviour, which would then allow one to test the hypothesis that some system is an agent by observing its behaviour.
(See \cite{Orseau2018}, which takes a Bayesian approach along these lines.)
Although Theorem \ref{thm.deterministicAreSpecifiable} is not a mathematically surprising result it does suggest a challenge that such an approach would need to overcome.

\subsection{Bounded rationality}
\label{sec.bounded}

Finally, we touch on the issue of bounded rationality.
A transducer may fail to be be optimal for a given teleo-environment $\environment$ while still being optimal for $\environment$ within some \emph{constrained class} of transducers, meaning simply a subset of all transducers.
These constrained classes of transducers are intended to model limitations in the cognitive capacity of a type of agent, such as a limited supply of memory, or of randomness.
For some kinds of limitation (i.e.\ some constrained classes of transducers), the value-laden Bellman property still holds.

In particular, the value-laden Bellman property holds for constrained classes
that are `closed under trajectory splicing', which is a somewhat technical property given in Definition \remoteReference{3.3.3} in the companion paper. For the interested reader, here is a brief informal description:
given transducers $\pi, \pi'$ and a sequence of inputs and outputs $(\mathbf{i},\mathbf{o})$, define a new transducer called ``$\pi$ spliced with $\pi'$ along $(\mathbf{i},\mathbf{o})$'', which behaves exactly like $\transducer$,
 except that if it happens to produce the sequence $\textbf{o}$ in response to the input sequence $\textbf{i}$
    starts behaving exactly as $\transducer'$ would have done in the same situation.
Then a constrained class of transducers $T$ is \emph{closed under trajectory splicing} if whenever $\pi,\pi\in T$, we also have that $\pi$ spliced by $\pi'$ along $(\mathbf{i},\mathbf{o})$ is in $T$, whenever $(\mathbf{i},\mathbf{o})$ is a possible sequence of inputs and outputs.

For agents that are optimal within such a class it is permissible to attribute normative-epistemic states that update by value-laden filtering, as described above.
However, this is not true for all constrained classes of interest.
We demonstrate this in Section \remoteReference{6.5} of the companion paper by applying our framework to a version of the `absent-minded driver,' \cite{PICCIONE19973} a classic problem from decision theory that involves making decisions without memory.
The optimal solution to the absent-minded driver problem is stochastic rather than deterministic, and as a transducer it is not possible to give it a mental state attribution that updates by value-laden filtering.

\section{Relation To Other Frameworks}
\label{sec.literature}

Our approach has some parallels with existing bodies of work in decision theory, cognitive neuroscience, machine learning and complex systems. In this section we discuss the relation of our approach to representation theorems in decision theory, the free-energy framework, inverse reinforcement learning, and computational mechanics.

\subsection{Representation Theorems For Expected Utility Theory}

It is helpful to compare our approach to classic `representation theorems' for
expected utility theory from the decision-theoretic literature. A representation
theorem (for expected utility) shows that, if a preference ordering $\prec$ over
`prospects'~$O$ has certain formal properties, then it is equivalent to a ranking
over prospects according to their expected value. As explained in \cite{sep-decision-theory} there are three well-known
representation theorems for expected utility: due to von Neumann / Morgenstern
(VNM) \cite{VNM1944}, Savage \cite{Savage1954}, and Bolker (used in Bolker-Jeffrey
decision theory). In VNM the prospects (the objects of the preference order) are `lotteries'; in Savage they are `acts'; and in Bolker they are `propositions'.

In some loose sense, we provide, in this paper, a representation result for deterministic behavioural propensities: we observe that they can be represented as the uniquely optimal solutions to control problems. However, the domain
in which we are working is very different from that of decision theory. We suppose that the object to be represented is the behaviour of a system; unlike VNM, Savage or Bolker/Jeffrey, we do not assume that information about an agent's preferences are directly available to the theorist. This is for reasons similar to those described in \cite{Hausman2000}: behaviour (the agent's `choices') does not reveal an agent's preferences unless we know the agent's beliefs, and beliefs in turn cannot be inferred from behaviour without information about the agent's evaluative judgements (which preferences are supposed to capture).

As discussed by \cite{Zynda2000}, the classic `representations' of preference orderings are non-unique, even given  a fixed form of representation.
Ours is no exception: there may be, and in general are, many different teleo-environments for which a given transducer is uniquely optimal.
However, in the classical representation frameworks, the representations are unique up to a class of fairly simple transformations.
This does not seem to be the case for our framework; there may be many teleo-environments that specify a given transducer, and these can be of very different character.

Another related family of classical results are the complete class theorems, for example \cite{wald1947essentially,Brown1981complete}.
These have the general character that, under differing assumptions, a decision rule with a property called being non-dominated also has the property of being `Bayes-optimal', meaning roughly that it can be thought of as doing loss minimisation with respect to some prior over the `states of nature.'
In our framework there is no obvious analog of the states of nature, so our results are of a different character; nevertheless it might be interesting to explore relationships in future work.

\subsection{The free energy principle}

An important contemporary area of research that touches on assigning mental states to systems is the free energy principle \cite{Friston2006brain,friston2019free}, also known as Bayesian mechanics \cite{Ramstead-2023}.
This is a broad body of theory that largely seeks applications in neuroscience.
However, it also applies to much more general physical systems.
Its proponents have talked about assigning mental states to non-biological physical systems since the earliest papers, e.g.\ \cite{Friston2006brain}; this is made very explicit in \cite{friston2019free}.
This body of work tends to use a different formal approach from the current work, focused on stochastic differential equations rather than stochastic processes over discrete alphabets.
Aside from technical differences, there are a few more foundational differences.
The free energy principle is concerned with internal states rather than only externally observable behaviour.
Perhaps more significantly, work on the free energy principle tends to seek a single, privileged belief state for a given system, which map states of the system to probability measures that relate to an `{actual}' or typical environment in which the system is found.
In the current paper we are focused on all possible assignments of belief states, including those in which the system might believe it inhabits a very different environment from the one it actually does;
for this reason, we do not consider the `true' environment at all in the current work.

\subsection{Inverse Reinforcement Learning}

Our work has some obvious similarities to an idea in machine learning known as inverse reinforcement learning.
Inverse reinforcement learning is the problem of inferring the reward function of an agent, given its
policy or observed behavior \cite{Arora2021}. Generally,  solutions to these problems assume that information about the environment is given, either
explicitly in the form of a transition function, or implicitly in the form of
trajectories. There are some partial exceptions to this, which attempt to infer
the reward function with much more limited knowledge of the
environment \cite{Choi2011, Klein2013}, but none that attempt to infer the whole
mental state of the agent based purely on its policy, which would be the true
equivalent of our approach. This is not surprising, as inverse reinforcement learning is famously an
already underspecified problem suffering from high computation costs, so the equivalent of our approach would likely be computationally infeasible.
We don't have such a constraint because we only aim to consider the abstract formal properties of such systems, rather than computing practical examples.

This concludes our review of previous work; in the next section we begin the main body of the paper.

\subsection{Computational Mechanics}

Our work is related to computational mechanics \cite{crutchfield2016complexity}, in that we are dealing with stochastic processes in discrete time, and our transducers have something somewhat in common with the causal states (i.e. the state of an $\varepsilon$-machine).
Our processes have inputs and outputs, but this can also be handled within computational mechanics \cite{barnett2015-epsilon}.
A bigger difference is that we do not consider stationarity at all, which makes the formalism quite different.
Our formalism doesn't make a distinction between stationary and nonstationary processes, because we  consider sequences that are only infinite in the future time direction, rather than the future and the past as in computational mechanics.
In addition, our primary focus is on the intentional stance, which hasn't previously been considered in the computational mechanics literature as far as we know.
We discuss the technical relationships to computational mechanics in more detail in the companion work.

\section{Limitations and future work}
\label{sec.limitations}

In this paper we outlined a formalism for thinking about several ways of
describing systems as potential agents and have proven some key results.
In this final section we discuss some of the limitations and possible extensions to the work, before concluding.

One can envisage several ways in which the work could be extended at the mathematical level. These include broadening or generalising the notion of norm (for example, in order to consider reward maximisation rather than success probability maximisation, or to consider multiple goals), as well as extensions to continuous time and continuous-valued inputs and outputs. Of particular interest would be considerations of multiple interacting systems, rather than only a single agent and its environment.
One could also consider the relationship between the \emph{agent's} beliefs about its environment (as attributed to it by the theorist), which we consider here, and the \emph{theorist's} beliefs about the same environment, which might be different.

In terms of limitations, one crucial missing component of this work is Dennett's idea that for some systems --- the co-called \emph{intentional systems} --- the intentional stance is much more productive than for others, and in particular more productive than the physical stance.
Our current framework only says which teleo-environments can be consistently attributed as beliefs; it doesn't currently offer any way to say which belief attributions are `better,' or more natural or pragmatic for a given agent's behaviour.
After all, in Theorem~\ref{thm.deterministicAreSpecifiable} we showed that one can always attribute a set of beliefs that say, in effect, that the agent believes it must act in exactly the way it does act, otherwise it will fail its goal; but we suspect it will rarely if ever be fruitful to explain a real agent's behaviour in terms of such a belief.
We suspect the answer might lie in bounded rationality, which we have only made a start on considering here.

Another, related, limitation of our approach is that we do not model so-called \emph{hyperintensional} dimensions of goals or beliefs.
For instance, it is possible to believe that there are seven biscuits in a jar without believing that the square of the number of biscuits in the jar is 49, although these propositions both refer to the same state of the world.
Roughly speaking, in our framework, when we attribute beliefs (i.e.\ a teleo-environment) to an agent, we are making a claim that the agent behaves consistently with those beliefs, but we are not making any claim that the agent \emph{thinks} or \emph{reasons about} those beliefs.
The things we call beliefs and goals in humans come equipped with a particular \emph{form} which is under-determined by the facts that they map to.
This notion of representational form, guise, or mode of presentation is related to the concept of information processing, which seems important if one wants to account for patterns of delays and errors in an agent's behaviour.
We suspect that addressing this issue would require significant additional formal machinery.

To conclude, we have given a simple framework in which some version of the intentional stance makes sense.
In particular we have shown that, within this framework, an agent that is optimal for some goal (i.e.\ teleo-environment) will in future time steps be optimal for an updated version of the same goal, where the updating is given by `value-laden (Bayesian) filtering.'
We have also shown that an agent is uniquely optimal for some teleo-environment if and only if it is deterministic.
However, there is much work to be done in extending the framework and in relating it more clearly to Dennett's ideas.

\section*{ACKNOWLEDGEMENT}
\small
This article was produced with financial and technical support from Principles of Intelligent Behaviour in Biological and Social Systems (PIBBSS), and Simon McGregor's work was supported with a scholarship grant from the Alignment of Complex Systems Research Group (ACS) at Charles University in Prague.
Nathaniel Virgo's work on this publication was made possible through the support of Grant 62229 from the John Templeton Foundation. The opinions expressed in this publication are those of the author(s) and do not necessarily reflect the views of the John Templeton Foundation.
The authors would also like to thank Martin Biehl for feedback on the manuscript.

\bibliography{notes}

\end{document}